\documentclass[preprint,12pt]{elsarticle}

\newtheorem{theorem}{Theorem}[section]
\newtheorem{lemma}{Lemma}[section]
\newtheorem{corollary}{Corollary}[section]

\newtheorem{remark}{Remark}[section]
\newcommand{\ignore}[1]{}{}

\usepackage{amssymb}
\usepackage{amsmath}

\usepackage{color}
\usepackage{soul}
\usepackage[dvipsnames]{xcolor}

\usepackage[colorlinks=true, linkcolor=blue, citecolor=blue]{hyperref}

\def\1{{{\mbox{${\rm{1\negthinspace\negthinspace I}}$}}}}

\renewcommand{\theequation}{\arabic{equation}}

\newcommand\beq{\begin{equation}}
\newcommand\eeq{\end{equation}}
\renewcommand{\theequation}{\thesection.\arabic{equation}}

\usepackage{ifthen}
\usepackage{xkeyval}
\usepackage{todonotes}
\begin{document}

\begin{frontmatter}

\title{Sharp moderate and large deviations for sample quantiles  }
\author{  Xiequan Fan  \ \      \ \      }
\address{School of Mathematics and Statistics, Northeastern University at Qinhuangdao, Qinhuangdao, China. }

\begin{abstract}
In this article, we discuss the sharp moderate and large deviations between
the quantiles  of population and the quantiles of samples.
 Cram\'{e}r type moderate deviations and Bahadur-Rao type  large deviations
are established with some mild conditions. The results refine the moderate and large deviation principles
of Xu and Miao [Filomat  2011; 25(2): 197-206].
\end{abstract}

\begin{keyword} sample quantiles, Cram\'{e}r moderate deviations, large deviation, sharp large deviations
\vspace{0.3cm}
\MSC primary 60E15; 60F10; secondary  62G30
\end{keyword}

\end{frontmatter}




\section{Introduction}\label{sec1}

Recall the definitions of the quantile of population  and the quantile of sample. Assume that
$\{X_i, 1\leq i \leq n  \}$ is a sample of a population $X$  with a common distribution function $F(x)$.
 For $0 < p < 1$,   denote
 $$x_p := F^{-1}(p)= \inf\{ x : F(x) \geq p \} $$
the $p$-th quantile of $F(x).$ An estimator of  $x_p$ is given by the sample $p$-th quantile defined as
follows:
$$x_{n,p} := F_n^{-1}(p)=  \inf\{x : F_n(x) \geq  p \},$$
where $$F_n(x) = \frac{1}{n} \sum_{i=1}^{n} \mathbf{1}(X_i \leq  x), \ \ \  x \in \mathbf{R},$$
is the empirical distribution function of the sample
$\{X_i, 1\leq i \leq n  \}$. Here,  $\mathbf{1}(A)$ denotes the indicator function of set $A.$

There are a number of literatures to study the asymptotic properties for sample quantiles.
If $x_{n,p}$ is the unique solution of the equation $F(x-)\leq p \leq F(x),$
 then  $x_{n,p} \rightarrow x_p$ a.e.;  Assume that $F(x)$ have a continuous
density function $f(x)$ in a neighborhood of $x_p$ such that  $f(x_p) >0$.
Then
 $  \frac{\sqrt{n} f(x_p) (x_{n,p} -x_p)}{\sqrt{p(1-p)} }  $ converges  to
 the standard normal random variable in distribution, see Serfling \cite{S80}.
Suppose that $F(x)$ is twice differentiable at  $x_p$, then Bahadur \cite{B66} proved
that
$$
x_{n,p}= x_p + \frac{p - F_n(x_p)}{f(x_p)} +R_n,\ \ \ a.e.,$$
where $R_n=O(n^{-3/4} (\ln n)^{3/4}) $ a.e., $n \rightarrow \infty.$
Xu and Miao \cite{XM06} (see also Miao, Chen  and Xu \cite{M11} for order statistics) obtained  the following moderate deviation principles for $x_{n,p} -x_p$. For
any sequence of real numbers  $\{a_n\}_{n\geq1}$   satisfying $a_n \rightarrow \infty$ and $a_n/\sqrt{n}\rightarrow 0$
as $n\rightarrow \infty$, it holds for any $r>0,$
\begin{eqnarray}\label{120ds1}
  \lim_{n\rightarrow \infty}\frac{1}{a_n^2}\ln \mathbf{P}\bigg(\frac{  \sqrt{n}}{a_n } |x_{n,p} -x_p|  \geq t \bigg) = -  \frac{  f(x_p)^2 r^2 }{ 2 p(1-p) } .
\end{eqnarray}
They also obtained the following  large deviation principles for $x_{n,p} -x_p$: the following two equality  hold  for any $t>0,$
 \begin{eqnarray}\label{120ds2}
  \lim_{n\rightarrow \infty}\frac{1}{n}\ln \mathbf{P}\bigg(   x_{n,p} -x_p   \geq t \bigg) = -  \inf_{y \geq 1-p}  \Lambda_+^*(y)
\end{eqnarray}
and
 \begin{eqnarray}\label{120ds3}
  \lim_{n\rightarrow \infty}\frac{1}{n}\ln \mathbf{P}\bigg(   x_{n,p} -x_p   \leq - t \bigg) = -    \inf_{y \geq  p}  \Lambda_-^*(y) ,
\end{eqnarray}
where
 \begin{eqnarray*}
  \Lambda_+^*(y) = y \ln  \frac{y}{F(x_p + t)}  +(1-y)\ln \frac{1-y}{1-F(x_p + t)}
\end{eqnarray*}
and
 \begin{eqnarray*}
 \Lambda_-^*(y) =y \ln  \frac{y}{F(x_p - t)}  +(1-y)\ln \frac{1-y}{1-F(x_p - t)} .
\end{eqnarray*}
In this paper, we are interested in sharp moderate and large deviations between
the quantiles of population and the quantiles of samples.   More precisely,
we establish Cram\'{e}r type moderate deviations and Bahadur-Rao type  large deviations
for $x_{n,p} -x_p$. Our results refine the moderate and large deviation principle  results  (\ref{120ds1}) - (\ref{120ds3}).

The paper is organized as follows. Our  results are stated and discussed in Section \ref{sec2}.
The proofs of theorems are deferred to Section \ref{sec3}.


\section{Main results}  \label{sec2}

For  brevity, denote
$$ R_n(x,p)= \frac{\sqrt{n} f(x_p) (x_{n,p} -x_p)}{\sqrt{p(1-p)} }, \ \ \ \ p \in (0, 1)   . $$
Denote   $\Phi(x)=\frac{1}{\sqrt{2\pi}}\int_{-\infty}^{x}\exp\{-t^2/2\}dt$   the standard normal distribution function.
The following theorem gives a Cram\'{e}r type  moderate deviation for sample quantiles.
\begin{theorem}\label{th00}
Let $f(x)$ be the density function of $X$ and let $p \in (0, 1).$
If $f'(x)$ is bounded in a neighbourhood of $x=x_p$ and  $f(x_p)> 0,$ then
 it holds
\begin{eqnarray}
  \ln \frac{\displaystyle \mathbf{P}\big( \pm  R_n(x,p)    \geq t   \big)}{1- \Phi(t)} = O\bigg( \frac{1+t^3}{\sqrt{n}}  \bigg)
\end{eqnarray}
uniformly for $0 \leq t =o( \sqrt{n} )  $
as $n\rightarrow \infty.$
\end{theorem}

\begin{remark}\label{remark2.1}
Let us comment on the result  of Theorem  \ref{th00}.
\begin{enumerate}
\item Assume that $f'(x)$ is uniformly bounded   on $\mathbf{R} $   and  that $f(x )$ is positive for all $x \in \mathbf{R}$.
When $p$ is replaced by $p_n$ which may depend  on $n,$ by inspecting the proof of Theorem  \ref{th00}, the following  equality
\begin{eqnarray}
  \ln \frac{\displaystyle \mathbf{P}\big( \pm  R_n(x,p_n)    \geq t   \big)}{1- \Phi(t)} = O\bigg( \frac{1+t^3}{\sqrt{n p_n(1-p_n)}\ }  \bigg)
\end{eqnarray}
holds uniformly for $0 \leq t =o( \sqrt{n  p_n(1-p_n)} \, )$ as $n\rightarrow \infty.$ Clearly, if $p_n(1-p_n) \rightarrow 0$ as $n\rightarrow \infty$, then the last range tends to smaller than $0 \leq t =o( \sqrt{n  } \, )$.

\item By an argument similar to the proof of Corollary 3 in \cite{FHM20},  the following moderate deviation principle (MDP)  result  is a  consequence of Theorem \ref{th00}.
Let   $\{a_n\}_{n\geq1}$ be  a  sequence of real numbers satisfying $a_n \rightarrow \infty$ and $a_n/\sqrt{n}\rightarrow 0$
as $n\rightarrow \infty$.  Then  for each Borel set $B$,
\begin{eqnarray}
- \inf_{x \in B^o}\frac{x^2}{2} &\leq & \liminf_{n\rightarrow \infty}\frac{1}{a_n^2}\ln \mathbf{P}\bigg(\frac{  R_n(x,p)}{a_n }  \in B \bigg) \nonumber \\
 &\leq& \limsup_{n\rightarrow \infty}\frac{1}{a_n^2}\ln \mathbf{P}\bigg(\frac{ R_n(x,p)}{a_n }    \in B \bigg) \leq  - \inf_{x \in \overline{B}}\frac{x^2}{2} ,   \label{MDP}
\end{eqnarray}
where $B^o$ and $\overline{B}$ denote the interior and the closure of $B$, respectively.
When $B$ is $[t, \infty)$ or $(-\infty, t]$ for some $t>0,$ the MDP result (\ref{MDP}) has been
established by Xu and Miao \cite{XM06} (cf.\, equality (\ref{120ds1})).  Notice that, in Xu and Miao \cite{XM06}, MDP result  holds without the assumption that
 $f'(x)$ is bounded in a neighbourhood of $x=x_p$.
\end{enumerate}
\end{remark}

Using the inequality $|e^x -1 | \leq   e^{ c } |x|$ valid for $|x| \leq c, $ from Theorem  \ref{th00},
we obtain  the following result about  the relative errors of normal approximations.
\begin{corollary}\label{corollary01}
Assume that the conditions of Theorem \ref{th00} are satisfied. Then  it holds
\begin{eqnarray}\label{ghkl5}
  \frac{\displaystyle \mathbf{P}\big( \pm  R_n(x,p)    \geq t \big)}{1- \Phi(t)} =1+  O \bigg( \frac{1+t^3}{\sqrt{n}} \bigg )
\end{eqnarray}
uniformly for $0 \leq t =O( n^{1/6})  $ as $n\rightarrow \infty,$  which  implies that
\begin{eqnarray}\label{ghkl6}
\frac{\displaystyle \mathbf{P}\big( \pm  R_n(x,p)    \geq t \big)}{1- \Phi(t)}=1+o(1)
\end{eqnarray}
uniformly for  $0\leq x = o ( n^{1/6} )$.
\end{corollary}

From (\ref{ghkl5}) in Corollary \ref{corollary01}, by an argument similar to the proof of Corollary 2.2 in Fan \emph{et al.}\,\cite{FGLS20}, we can obtain the following Berry-Esseen bound:
\begin{eqnarray}
\sup_{t \in \mathbf{R}} \Big| \mathbf{P}\big( \pm  R_n(x,p)  \leq t   \big) -  \Phi(t) \Big| = O\Big( \frac{1}{\sqrt{n}}  \Big),\ \ \ n\rightarrow \infty.
\end{eqnarray}
The last convergence rate coincides with the classical result established by Reiss \cite{R74}.

Theorem \ref{th00} is devoted to the moderate deviations. For sharp large deviations,   we have the following   Bahadur-Rao type large deviation expansions.
\begin{theorem}\label{th02}
For any $t \geq 0$, it holds
\begin{eqnarray*}
 \mathbf{P}\big(  x_{n,p} -x_p    \geq  t   \big) =\frac{1}{\tau_t^+ \, \sigma _p   \sqrt{2\pi n  }\ } e^{-n \Lambda^+(t)}\Big[1+o(1) \Big], \ \ \ n \rightarrow \infty,
\end{eqnarray*}
where
\begin{eqnarray*}
&&\tau_t^+ = \ln \frac{ F(x_p + t) (1-p)}{p(1-F(x_p + t))}  ,\ \ \   \sigma_p =  \sqrt{p(1-p) } \ \ \ \ \textrm{and} \  \\
 \nonumber \\
&& \Lambda^+(t)=p \ln  \frac{p}{F(x_p + t)}  +(1-p)\ln \frac{1-p}{1-F(x_p + t)} .
\end{eqnarray*}
Similarly, it also holds for any $t \geq 0$,
\begin{eqnarray*}
 \mathbf{P}\big(  x_{n,p} -x_p    \leq  -  t   \big)  =\frac{1}{ \tau_t^-\,\sigma _p   \sqrt{2\pi n }\ }e^{-n \Lambda^-(t)}\Big[1+o(1) \Big],\ \ \ n \rightarrow \infty,
\end{eqnarray*}
where
\begin{eqnarray*} \label{lamda2}
&&\tau_t^- = \ln \frac{p(1-F(x_p - t))}{F(x_p - t)(1-p)}  \ \ \  \ \   \textrm{and}   \\
 \nonumber \\
&& \Lambda^-(t)=  p \ln  \frac{p}{F(x_p - t)}  +(1-p)\ln \frac{1-p}{1-F(x_p - t)}.
\end{eqnarray*}
\end{theorem}

Denote  $c$ a finite and positive constant  which does not depend on $n$.    For  two sequences of positive numbers $(a_n)_{n\geq 1}$ and $(b_n)_{n\geq 1}$,   write $a_n \asymp b_n$ if there exists a   $c $ such that ${a_n}/{c}\leq b_n\leq c\, a_n$ for all sufficiently large $n$.
By Theorem \ref{th02},  we have for any given constants $t > 0$ and $p \in (0, 1)$,
\begin{eqnarray}
 \mathbf{P}\big(  x_{n,p} -x_p    \geq  t   \big) \asymp    \frac{1}{\sqrt{n} } e^{-n \Lambda^+(t)} \ \ \
\textrm{and}
\ \ \
 \mathbf{P}\big(  x_{n,p} -x_p    \leq  -  t   \big) \asymp  \frac{1}{\sqrt{n}    } e^{-n \Lambda^-(t)}, \ \ \ n  \rightarrow \infty. \nonumber
\end{eqnarray}
In particular, from the last line, we recover the following large deviation principle (LDP) result of Xu and Miao \cite{XM06}:
for any $t>0,$ the following equalities hold
\begin{eqnarray}
  \lim_{n\rightarrow \infty}\frac{1}{n}\ln \mathbf{P}\big(  x_{n,p} -x_p    \geq  t   \big) = - \Lambda^+(t) \ \ \ \
\textrm{and} \ \ \ \
  \lim_{n\rightarrow \infty}\frac{1}{n}\ln \mathbf{P}\big(  x_{n,p} -x_p   \leq-  t   \big) =  - \Lambda^-(t), \nonumber
\end{eqnarray}
see equalities (\ref{120ds2}) - (\ref{120ds3}). Notice that $\Lambda^+(t)= \inf_{y \geq 1-p}  \Lambda_+^*(y)$ and $  \Lambda^-(t)=  \inf_{y \geq  p}  \Lambda_-^*(y) ,$ see Remark 1 in Xu and Miao \cite{XM06}.

\section{Proofs of Theorems} \label{sec3}
\setcounter{equation}{0}

\subsection{Proof  of Theorem \ref{th00}}
Let $(Y_{i})_{i\geq 1} $ be a sequence of i.i.d.\ and centered random variables.
Denote $\sigma^2=\mathbb{E}Y_{1}^2$ and $T_n=\sum_{i=1}^{n}Y_{i}.$
Cram\'{e}r \cite{Cramer38} has established the following asymptotic expansion on the  probabilities of moderate deviations for $T_n$.
\begin{lemma}\label{lem3.1}
Assume that $\mathbf{E} e^{ \lambda |Y_1|}  < \infty  $
for a constant $\lambda>0$. Then it holds
\begin{equation}
\ln \frac {\mathbb{P}( T_n> x\sigma\sqrt{n})} {1-\Phi(x)}  =  O \bigg( \frac{1+x^3}{\sqrt{n}}\bigg) \ \ \mbox{as} \ \ n \rightarrow \infty,
\label{cramer1}
\end{equation}
 uniformly for $0\leq x =  o(n^{1/2}).$
\end{lemma}
The Cram\'{e}r moderate deviations  have attracted a lot of interests.  We refer to  Petrov \cite{Petrov75}, Beknazaryan, Sang and Xiao \cite{BXY19}, Fan, Hu and Ma \cite{FHM20}, 
Fan, Hu and Xu \cite{FHX22}
for more such type results.

 We first give a proof for the case  $R_n(x,p)$.  For all $t \geq 0,$ it is easy to see that
\begin{eqnarray*}
  \mathbf{P}\big(  R_n(x,p)    \geq t \big) &=&  \mathbf{P}\bigg(  \frac{\sqrt{n} f(x_p) (x_{n,p} -x_p)}{\sqrt{p(1-p)} }   \geq t \bigg)  \\
  &=& \mathbf{P}\bigg( x_{n,p} -x_p     \geq  \frac{t\sqrt{p(1-p)} }{\sqrt{n} f(x_p)} \bigg).
\end{eqnarray*}
Write $x_{n,p,t} = x_p +  \frac{t\sqrt{p(1-p)} }{\sqrt{n} f(x_p)}.$ Then, by the definition of $x_{n,p,t}$, we get for all $t \geq 0,$
\begin{eqnarray*}
  \mathbf{P}\big(  R_n(x,p)    \geq t \big)  &=& \mathbf{P}\big( x_{n,p}     \geq  x_{n,p,t} \big) \ =\ \mathbf{P}\big( p    \geq F_n( x_{n,p,t}) \big) \\
  &=& \mathbf{P}\Big(  n p    \geq \sum_{i=1}^{n}  \mathbf{1}(X_i \leq x_{n,p,t}  ) \Big) \\
    &=& \mathbf{P}\bigg(  n p-nF(x_{n,p,t})    \geq \sum_{i=1}^{n} \big( \mathbf{1}(X_i \leq x_{n,p,t}   ) - F(x_{n,p,t}) \big ) \bigg) \\
     &=& \mathbf{P}\Bigg( \sum_{i=1}^{n}  \frac{\mathbf{1}(X_i \leq x_{n,p,t}   ) - F(t_{n,p}(x))}{  \sqrt{ n F(x_{n,p,t})(1-F(x_{n,p,t}))}  }  \leq
     \frac{\sqrt{n} ( p- F(x_{n,p,t}))}{\sqrt{  F(x_{n,p,t})(1-F(x_{n,p,t}))} } \Bigg).
\end{eqnarray*}
Recall that $ F'(x_{p})= f(x_{p})$ and that $F''(x)=f'(x)$ is bounded in a neighbourhood of $x=x_p$. Thus, it holds uniformly for  $0 \leq t =o( \sqrt{n}),$
\begin{eqnarray*}
p- F(x_{n,p,t})   &=& F(x_{p})- F(x_{n,p,t})  \\
&=& -  \frac{t }{\sqrt{n}  }\sqrt{p(1-p)} + O\bigg(\frac{t^2 }{n  }   \bigg).
\end{eqnarray*}
From the last line, we deduce that uniformly for  $0 \leq t =o( \sqrt{n}),$
\begin{eqnarray}\label{gsdc2}
F(x_{n,p,t}) \big(1- F(x_{n,p,t}) \big)
&=&p(1-p) + \frac{(1-2p)t }{\sqrt{n p(1-p)}  }  p(1-p) + O\bigg(\frac{t^2 }{n  }   \bigg).
\end{eqnarray}
Hence, we have uniformly for  $0 \leq t =o( \sqrt{n}),$
\begin{eqnarray*}
\frac{\sqrt{n} ( p- F(x_{n,p,t}))}{\sqrt{  F(x_{n,p,t})(1-F(x_{n,p,t}))} }
&=&-t + O\bigg(\frac{t^2 }{\sqrt{np(1-p)}  }   \bigg).
\end{eqnarray*}
Therefore, we deduce that uniformly for  $0 \leq t =o( \sqrt{n}),$
\begin{eqnarray*}
  \mathbf{P}\big(  R_n(x,p)    \geq t \big)
      =  \mathbf{P}\Bigg( \sum_{i=1}^{n}  \frac{\mathbf{1}(X_i \leq x_{n,p,t}   ) - F(t_{n,p}(x))}{  \sqrt{ n F(x_{n,p,t})(1-F(x_{n,p,t}))}  }  \leq
     -t + O\bigg(\frac{t^2 }{\sqrt{np(1-p)}  }   \bigg) \Bigg).
\end{eqnarray*}
Denote $Z_i= \mathbf{1}(X_i \leq x_{n,p,t}   ) - F(t_{n,p}(x)), 1 \leq i \leq n.$
  Notice that $ (Z_i  )_{1 \leq i \leq n}$ are i.i.d.\ and centered random variables,
   and satisfy that for all $1\leq i\leq n,$
$$ | Z_i| \leq 1 \ \ \ \ \ \textrm{and} \ \ \ \ \   \textrm{Var}   ( Z_i  ) =  F(x_{n,p,t})\big(1-F(x_{n,p,t}) \big).$$
The last line and (\ref{gsdc2}) implies that
\begin{eqnarray} \label{fsddf}
\sum_{i=1}^n \textrm{Var}   ( Z_i  )=nF(x_{n,p,t}) \big(1- F(x_{n,p,t}) \big)
 = &np(1-p) +   O\big( \sqrt{n}     \big).
\end{eqnarray}
By Lemma \ref{lem3.1} and (\ref{fsddf}), we obtain that it holds
\begin{equation}
\ln \frac { \mathbf{P}\big(  R_n(x,p)    \geq t \big) } {1-\Phi(t)}  \ =\  O \bigg( \frac{1+ t^3}{\sqrt{n}}\bigg) \ \ \mbox{as} \ \ n \rightarrow \infty,
\end{equation}
 uniformly for $0\leq t =  o(n^{1/2}), $ which gives the  desired inequality for $R_n(x,p)$. For $-  R_n(x,p)$, the desired inequality follows by a similar argument.    \hfill\qed

%
%
%

\subsection{Proof  of Theorem \ref{th02}}
Recall that $(Y_{i})_{i\geq 1} $ is a sequence of i.i.d.\,and centered random variables and $T_n=\sum_{i=1}^{n}Y_{i}.$  Assume that $\mathbf{E} e^{ \lambda |Y_1|}  < \infty  $
for a constant $\lambda>0$.
Denote $$\Lambda ^\ast(x)=\sup_{\lambda\geq 0}\{\lambda x -  \log \mathbf{E} e^{\lambda Y_1} \}$$   the Fenchel-Legendre transform  of the  cumulant  function of $Y_1$. The function $\Lambda^\ast(x) $ is   known as the good rate function in  LDP theory, see   Dembo and Zeitouni \cite{D98}.
Bahadur and Rao \cite{BR60} have established the following sharp large deviations.

\begin{lemma}
\label{thend}
Assume that $\mathbf{E} e^{ \lambda |Y_1|}  < \infty  $
for a constant $\lambda>0$. For any $y>0,$ let $ \tau_y  $ and $\sigma_y $ be the positive solutions of  the following equations:  $$h'(\tau_y)=0 \ \ \ \ \textrm{and}\ \ \ \ \sigma_y=\sqrt{ -h''(\tau_y) } ,$$
where $h(\tau)=\tau y-\log  \mathbf{E}  e^{\tau Y_1}.$ Then for a given positive constant $y$, it holds
\begin{eqnarray}
\mathbf{P}\bigg( \frac{T_{n}}{n} > y \bigg)= \frac{  e^{- n \, \Lambda^\ast(y)} }{  \sigma _y  \, \tau_y   \sqrt{2\pi n }  } \Big[ 1+ o(1) \Big], \ \ \ \ \ n\rightarrow \infty.
\end{eqnarray}
\end{lemma}
Such  type large deviations have attracted a lot of attentions.  We refer to Bercu and  Rouault \cite{BR06}, Joutard \cite{J06}, Fan, Grama and Liu \cite{FGL15}  for more
such type results.

We are in position to prove Theorem \ref{th02}.
For all $t \geq 0,$ it holds
\begin{eqnarray*}
  \mathbf{P}\big(  x_{n,p} -x_p    \geq t   \big)  &=& \mathbf{P}\big( x_{n,p}     \geq  x_p + t \big) \ =\ \mathbf{P}\big( p    \geq F_n(x_p + t ) \big) \\
  &=& \mathbf{P}\Big(  n (1-p)   \leq \sum_{i=1}^{n}  \mathbf{1}(X_i > x_p + t ) \Big) \\
    &=& \mathbf{P}\bigg(     n ( F(x_p + t)-p)  \leq \sum_{i=1}^{n} U_i \bigg) ,
\end{eqnarray*}
where
$$U_i= \mathbf{1}(X_i > x_p + t   ) -1+ F(x_p + t),\ \ \ \ \ \ i=1,...,n.$$
Notice that $(U_i)_{1\leq i \leq n}$ are i.i.d.\ and centered random variables with $|U_i|\leq 1.$
By some simple calculations,  it is easy to see that
\begin{eqnarray*}
  \mathbf{E} e^{\lambda  U_1}
  &=&  e^{\lambda  F(x_p + t)}\big(1-F(x_p + t) \big)+  e^{\lambda ( F(x_p + t)-1)} F(x_p + t)    \\
  &=&\exp\bigg\{   \lambda  F(x_p + t)+  \ln  \Big(  1-F(x_p + t) +  e^{-\lambda } F(x_p + t)  \Big)    \bigg\}.
\end{eqnarray*}
Therefore, we have
\begin{eqnarray*}
 \Lambda ^\ast( F(x_p + t)-p) &=& \sup_{\lambda\geq 0} \Big\{\lambda ( F(x_p + t)-p) -  \log \mathbf{E} e^{\lambda U_1} \Big\} \\
 &=& \sup_{\lambda\geq 0}\Big\{- \lambda p   - \ln  \Big(  1-F(x_p + t) +  e^{-\lambda } F(x_p + t)  \Big)  \Big\} \\
  &= &    \Lambda^+(t)    ,
\end{eqnarray*}
where
$$ \Lambda^+(t)=p \ln  \frac{p}{F(x_p + t)}  +(1-p)\ln \frac{1-p}{1-F(x_p + t)}. $$
Denote $$h_1(\tau)=\tau ( F(x_p + t)-p) -\log  \mathbf{E}  e^{\tau U_1}.$$
Let $ \tau_t^+  $ and $\sigma_t^+ $ be the positive solutions of  the following equations:  $$h_1'(\tau_t^+)=0 \ \ \ \ \textrm{and}\ \ \ \ \sigma_t^+=\sqrt{ -h_1''(\tau_t^+) } .$$
Then we have
 $$ \tau_t^+ = \ln \frac{ F(x_p + t) (1-p)}{p(1-F(x_p + t))}     \ \  \ \ \textrm{and} \ \ \  \ \ \sigma_t^+=\sigma_p=  \sqrt{p(1-p) }  .$$
Applying Lemma \ref{thend}  to $ (U_{i})_{1 \leq i \leq n} $, with $$ Y_i=U_i, \ \ \  y=F(x_p + t)-p, \ \ \ \tau_t  =  \tau_t^+  \ \ \  \textrm{and} \ \ \   \sigma_t= \sigma_p , $$
we get for any $ t >0$,
\begin{eqnarray*}
 \mathbf{P}\big(  x_{n,p} -x_p    \geq t   \big)    =  \frac{1}{ \sigma _p \, \tau_t^+   \sqrt{2\pi n }\ }  \exp\Big\{ -n \Lambda^+(t)     \Big\}\Big[ 1+ o(1) \Big], \ \ \ n \rightarrow \infty,
\end{eqnarray*}
which gives the first desired equality.

For any $t > 0,$ it holds
\begin{eqnarray*}
  \mathbf{P}\big(  x_{n,p} -x_p \leq - t   \big)  &=& \mathbf{P}\big( x_{n,p} \leq  x_p - t \big) \ =\ \mathbf{P}\big( p    \leq F_n(x_p - t ) \big) \\
  &=& \mathbf{P}\Big(  n p   \leq \sum_{i=1}^{n}  \mathbf{1}(X_i \leq x_p - t ) \Big) \\
    &=& \mathbf{P}\bigg(     n (p- F(x_p -t))  \leq \sum_{i=1}^{n} V_i \bigg) ,
\end{eqnarray*}
where
$$V_i= \mathbf{1}(X_i \leq x_p - t   ) -  F(x_p - t),\ \ \ \ \ \ i=1,...,n.$$
Notice that $(V_i)_{1\leq i \leq n}$ are i.i.d.\ and centered random variables with $|V_i|\leq 1.$
By some calculations,  it is easy to see that
\begin{eqnarray*}
  \mathbf{E} e^{\lambda  V_1}
  &=& e^{\lambda (1-  F(x_p- t)) } F(x_p -t)  +  e^{-\lambda   F(x_p - t) } (1-F(x_p - t))     \\
  &=&\exp\bigg\{  -  \lambda  F(x_p - t)+  \ln  \Big(  1-F(x_p - t) +  e^{\lambda } F(x_p - t)  \Big)     \bigg\}.
\end{eqnarray*}
Therefore, we have
\begin{eqnarray*}
 \Lambda ^\ast( p- F(x_p -t)) &=& \sup_{\lambda\geq 0} \Big\{\lambda ( p- F(x_p -t)) -  \log \mathbf{E} e^{\lambda V_1} \Big\} \\
 &=& \sup_{\lambda\geq 0}\Big\{ \lambda p   - \ln  \Big(  1-F(x_p + t) +  e^{-\lambda } F(x_p + t)  \Big)  \Big\} \\
  &= &    \Lambda^-(t)    ,
\end{eqnarray*}
where
$$ \Lambda^-(t)= p \ln  \frac{p}{F(x_p - t)}  +(1-p)\ln \frac{1-p}{1-F(x_p - t)}.  $$
Denote $$h_2(\tau)=\tau ( p- F(x_p -t)) -\log  \mathbf{E}  e^{\tau V_1}.$$
Let $ \tau_t^-  $ and $\sigma_t^- $ be the positive solutions of  the following equations:  $$h_2'(\tau_t^-)=0 \ \ \ \ \textrm{and}\ \ \ \ \sigma_t^-=\sqrt{ -h_2''(\tau_t^-) } .$$
Then we have
 $$ \tau_t^- =\ln \frac{p(1-F(x_p - t))}{F(x_p - t)(1-p)}    \ \  \ \ \textrm{and} \ \ \  \ \ \sigma_t^-=\sigma_p= \sqrt{p(1-p) }.$$
Applying Lemma \ref{thend}  to $ (V_{i})_{1 \leq i \leq n} $, with $$ Y_i=V_i, \ \ \  y=p- F(x_p -t), \ \ \ \tau_t  =  \tau_t^-  \ \ \  \textrm{and} \ \ \   \sigma_t= \sigma_p , $$
we get for any $ t >0$,
\begin{eqnarray*}
 \mathbf{P}\big(  x_{n,p} -x_p    \leq - t    \big)    =  \frac{1}{ \sigma _p \, \tau_t^-   \sqrt{2\pi n }\ }  \exp\Big\{ -n \Lambda^-(t)     \Big\}\Big[ 1+ o(1) \Big], \ \ \ n \rightarrow \infty,
\end{eqnarray*}
which gives the second desired equality.


\end{document}